\documentclass[12pt]{article}
\usepackage{euscript,amsmath, amssymb, amsfonts}
\usepackage{color}

\newcommand{\R}{{\cal R}}

\newcommand{\G}{{\mathbb C}[W(\R)]}

\newcommand{\be}{\begin{equation}}
\newcommand{\ee}{\end{equation}}
\newcommand{\bee}{\begin{eqnarray}}
\newcommand{\eee}{\end{eqnarray}}
\newcommand\nn{\nonumber \\}

\newcommand{\yy}{{\vec y}}

\newcommand{\vv}{{\vec v}}

\newcommand{\strr}{ {\mathrm str} }
\newcommand{\trr}{ {\mathrm tr} }

\newenvironment{proof}[1][Proof]{\noindent\textsf{#1.} }{\ \rule{0.5em}{0.5em}}

\newcounter{theorem}

\newcounter{lemma}

\newcounter{proposition}

\newcounter{definition}

\font\frtnfr=eufm10   scaled\magstep1 \font\twlfr=eufm10
\font\tenfr=eufm10  
\newfam\frfam
\textfont\frfam=\frtnfr \scriptfont\frfam=\twlfr
\scriptscriptfont\frfam=\tenfr

 \font\twlopen=msbm10
\font\tenopen=msbm10  
\newfam\openfam
\textfont\openfam=\tenopen \scriptfont\openfam=\twlopen
\scriptscriptfont\openfam=\tenopen

\font\frtnsf = cmss12 scaled\magstep1 \font\twlsf = cmss10
\font\tensf = cmss9  
\newfam\Scfam
\textfont\Scfam = \frtnsf \scriptfont\Scfam = \twlsf
\scriptscriptfont\Scfam = \tensf

\begin{document}

\renewcommand{\theequation}{\arabic{equation}}
\bibliographystyle{nphys}
%%%%%%%%%%%%%%%%%%%%%%%%%%%%%%%%%%%%%%%%%%%%%%%%%%%%%%%%%

%%%%%%%%%%%%%%%%%%%%%%%%%%%%%%%%%%%%%%%%%%%%%%%%%%%%%%%%%%%%%%%%%%%%%
%%%%%%%%%%%%%%      TITLE PAGE     %%%%%%%%%%%%%%%%%%%%%%%%%%%%%%%%%%
%%%%%%%%%%%%%%%%%%%%%%%%%%%%%%%%%%%%%%%%%%%%%%%%%%%%%%%%%%%%%%%%%%%%%

\sloppy

\title{      Klein operator and the Number of independent Traces and Supertraces
      on the Superalgebra of Observables of Rational Calogero Model\\
      based on the Root System
}

\author
 {
 S.E. Konstein\thanks{E-mail: konstein@lpi.ru}
               \thanks{ I.E.Tamm Department of
               Theoretical Physics, P. N. Lebedev Physical
               Institute, Leninsky Prospect 53,
               119991 Moscow, Russia.}
  and
 R. Stekolshchik \thanks{E-mail: r.stekol@gmail.com}
  }

\date{}

\maketitle
\thispagestyle{empty}

\begin{abstract}
In the Coxeter group $W(\R)$ generated by the root system $\R$, let
$ T(\R)$ be the number of conjugacy classes having no eigenvalue $+1$
and let $ S(\R)$ be the number of conjugacy classes having no
eigenvalue $-1$. The algebra $H_{W(\R)}$ of observables of the
rational Calogero model based on the root system $\R$ possesses
$ T(\R)$ independent traces; the same algebra, considered as an associative
superalgebra with respect to a certain natural parity, possesses
$ S(\R)$ even independent supertraces and no odd trace or supertrace. The numbers
$ T(\R)$ and $ S(\R)$ are determined for all irreducible root systems
(hence for all root systems). It is shown that $ T(\R)\le  S(\R)$, and
$ T(\R) =  S(\R)$ if and only if superalgebra $H_{W(\R)}$ contains a
Klein operator (or, equivalently, $W(\R) \ni -1$).
\end{abstract}

%\keywords{Trace; supertrace; Cherednik algebra; algebra of observables; Calogero model.}
%\ccode{2000 Mathematics Subject Classification: 17B80, 16W55}

%%%%%%%%%%%%%%%%%%%%%%%%%%%%%%%%%%%%%%%%%%%%%%%%%%%%%%

\section{Definitions and generalities}

\subsection{Traces}

Let ${\cal A}$ be an associative superalgebra with parity $\pi$. All
expressions of linear algebra are given for homogeneous elements
only and are supposed to be extended to inhomogeneous elements via
linearity.

%\begin{definition}\label{str}
A linear function $str$ on ${\cal A}$ is called a {\it supertrace} if
$$str(fg)=(-1)^{\pi(f)\pi(g)}str(gf) \ \mbox{ for all } f,g\in {\cal A}.$$
%\end{definition}

%\begin{definition}\label{tr}
A linear function $tr$ on ${\cal A}$ is called a {\it trace} if
$$tr(fg)=tr(gf) \ \mbox{ for all } f,g\in {\cal A}.$$
%\end{definition}

%\medskip

%\begin{definition}
 A linear function $L$ is {\it even} if $L(f)=0$
for any $f\in{\cal A}$ such that $\pi(f)=1$,
it is
{\it odd} if $L(f)=0$
for any $f\in{\cal A}$ such that $\pi(f)=0$.
%\end{definition}

%\medskip
Let ${\cal A}_1$ and
${\cal A}_2$ be associative superalgebras with parities $\pi_1$ and $\pi_2\,$,
respectively.
Define the tensor product
${\cal A}={\cal A}_1\otimes {\cal A}_2$
as a superalgebra with the product
$(a_1 \otimes a_2)(b_1 \otimes b_2)=(a_1 b_1) \otimes (a_2 b_2)$ (no sign factors in this formula)
and the
parity $\pi$ defined by the formula  $\pi(a\otimes b)=\pi_1(a)+\pi_2(b)$.

Let $T_i$ be a trace on  ${\cal A}_i$.
Clearly, the function $T$ such that $T(a\otimes b)=T_1(a)T_2(b)$ is a trace on ${\cal A}$.

Let $S_i$ be an {\it even} supertrace on  ${\cal A}_i$. Clearly,
the function $S$ such that $S(a\otimes b)=S_1(a)S_2(b)$ is an even
supertrace on ${\cal A}$.

\subsection{Klein operator}

{Let ${\cal A}$ be an associative superalgebra with parity $\pi$.
Following M.Vasiliev, see, e.g. \cite{V}, we say that an element
$K\in {\cal A}$ is a {\it Klein operator}\footnote{In honor of Oskar
Klein.}
 if $\pi(K)=0$,
$Kf=(-1)^{\pi(f)}fK$ for any $f\in {\cal A}$ and $K^2=1$. Every Klein operator
belongs to the {\it anticenter} of the superalgebra ${\cal A}$, see \cite{sosrus}, p.41.%
\footnote{Let ${\cal A}$ be an associative superalgebra with parity
$\pi$. Its {\it anticenter $AC({\cal A})$} is defined by the formula
$$
AC({\cal A})=\{
a\in {\cal A} \,\,|\,\, ax-(-1)^{\pi(x)(\pi(a)+1)}xa =0
 \ \text{ for any } x\in{\cal A }\}.
$$
}

Any Klein operator, if exists, establishes an isomorphism between the space of
even traces and the space of even supertraces on ${\cal A}$.
Namely, if $f\mapsto T(f)$ is an even trace, then $f\mapsto T(fK)$
is a supertrace, and if $f\mapsto S(f)$ is an even supertrace, then $f\mapsto S(fK)$ is a trace.

%%%%%%%%%%%%%%%%%%%%%%%%%%%%%%%%%%%%%%%%%%%%%%%%%%%%%%%%%%%%%%%%%%%

\subsection{Group algebra}

Let $V={\mathbb  R}^n$ and $G\subset End(V)$ be a finite group. The
{\it group algebra} ${\mathbb C}[G]$  of $G$ consists of all linear
combinations $\sum_{g\in W(\R)} \alpha_g \bar g$, where $\alpha_g
\in {\mathbb C}$. We distinguish $g$ considered as an element of the
group $G\subset End(V)$ from the same element $\bar g \in {\mathbb
C}[G]$ considered as an element of the group algebra. The addition
in ${\mathbb C}[G]$ is defined as follows:
$$
\sum_{g\in G} \alpha_g \bar g + \sum_{g\in G} \beta_g \bar g
= \sum_{g\in G} (\alpha_g + \beta_g) \bar g
$$
and the multiplication is defined by setting
$\overline {g_1\!}\,\, \overline {g_2\!} = \overline {g_1 g_2}$.

Note that the additions in ${\mathbb C}[G]$ and in $End (V)$ differ.
For example, if $I\in G$ is unity and the matrix $J=-I$ from $End(V)$
belongs to $G$, then $I+J=0$ in $End(V)$ while
$\overline I + \overline J \ne 0$ in ${\mathbb C}[G]$.

\subsection{Root systems}

%%%%%%%%%%%%%%%%%%%%%%%%%%%%%%%%%%%%%%%%%%%%%%%%%%%%%%%%%%%%%%%%%%%%%%%%%%%%%%%%%%%%%%%
%\begin{definition}\label{reflec}
Let $V={\mathbb  R}^N$ be endowed with a non-degenerate symmetric
bilinear form $(\cdot,\cdot)$ and the vectors $\vec a_i$ constitute
an orthonormal basis in $V$, i.e.
$$(\vec a_i,\, \vec a_j)=\delta_{ij}.$$
Let
$x^i$ be
the coordinates of ${\vec x} \in V$, i.e.  ${\vec x} =\vec a_i\,x^i$.
Then $({\vec x} ,\,\yy)=\sum_{i=1}^N x^i y^i$ for any ${\vec x} ,\,\yy \in V$.
The indices $i$ are raised and lowered by means of the forms $\delta_{ij}$
and $\delta^{ij}$.

For any nonzero $\vv \in V={\mathbb  R}^N$, define the {\it
reflections} $R_\vv$ as follows: \be\label{ref} R_\vv ({\vec x}
)={\vec x}  -2 \frac {({\vec x} ,\,\vv)} {(\vv,\,\vv)} \vv \qquad
\mbox{ for any }{\vec x}  \in V. \ee

The reflections (\ref{ref}) have the following properties
\be
%\label{prop}
R_\vv (\vv)=-\vv,\qquad R_\vv^2 =1,\qquad
({R}_\vv ({\vec x} ),\,\vec u)=({\vec x} ,\,{R}_\vv (\vec u))\quad
\mbox{ for any }\vv,\,{\vec x} ,\,\vec u\in V.\nonumber
\ee

%\begin{definition}\label{rootsys}
A finite set of vectors $\R\subset V$ is said to be a {\it root system}
if the following conditions hold:

i) $\R$ is ${R}_\vv$-invariant for any $\vv \in \R$,

ii) if $\vv_1,\vv_2\in \R$ are collinear, then either $\vv_1=\vv_2$ or $\vv_1=-\vv_2$.
%\end{definition}

%\medskip
Clearly, the group $W(\R)\subset O(N, {\mathbb  R})\subset End(V)$
generated by all reflections ${R}_\vv$ with $\vv \in \R$ is finite.

%%%%%%%%%%%%%%%%%%%%%%%%%%%

%\medskip

Let $V=V_1 \oplus V_2$, where $V_1\ne \{0\}$ and $V_2\ne \{0\}$
are orthogonal with respect to the form $(\cdot,\cdot)$,
and let a root system on $V$ have a decomposition:
$\R=\R_1 \bigcup \R_2$, where $\R_i \subset V_i$ for $i=1,2$.
Then each $\R_i\subset V_i$ is a root system.
We say in this case that $\R$ is {\it reducible}, and denote this fact
as $\R=\R_1 + \R_2$.
Note, that each $\R_i$ can be empty.
A root system which is not reducible is called {\it irreducible}.

If  $\R=\R_1 + \R_2$, then $W(\R)=W(\R_1)\times W(\R_2)$.

Any root system has a decomposition $\R=\R_1+\R_2+ \ldots +\R_n$,
where the $\R_j$ are irreducible root systems.

All irreducible root systems are listed in numerous literature
(see, e.g., \cite{NB}, \cite{JH}, \cite{OP}, \cite{carter}, \cite{carst}, \cite{h4}, \cite{raf}).
As it follows from the definition of a root system given above, we consider
both crystallographic ($A_n$, $B_n$, $C_n$, $D_n$, $E_6$, $E_7$, $E_8$, $F_4$, $G_2$)
and non-crystallographic ($H_3$, $H_4$, $I_2(n)$) root systems.

We consider also the empty root system,
assuming that it generates the trivial group consisting of the unity element only.

The definition of reducible root system implies that the empty root
system in ${\mathbb  R}^N$ is reducible for any $N>1$. The
{\it irreducible empty root system} --- we  denote it $A_0$ --- belongs to
${\mathbb  R}$.

%%%%%%%%%%%%%%%%%%%%%%%%%%%%%%%%%%%%%%%%%%%%%%%%%%%%%%%%%%%%%%%%%%%%%%%%%%%%%%%%%%%%%%%%%%%%%%%%%%%%%%%%%

\subsection{The superalgebra of observables}

Let $\R$ be a finite root system. Let $\eta$ be a set of constants
$\eta_\vv$ with $\vv\in\R$ such that $\eta_\vv=\eta_{\vec w}$ if
$R_\vv$ and $R_{\vec w}$ belong to one conjugacy class of $W(\R)$.

Let ${\cal H}^\alpha$, where $\alpha=0,1$\,, be two copies of $V$
with orthonormal bases $a_{\alpha\,i} \in {\cal H}^\alpha$, where
$i=1,\,...\,,\,N$.

 %\medskip

\begin{definition}%\label{HW}
%{\bf Definition.}
%\label{defpage}
The superalgebra $H_{W({\R})}(\eta)$ is an associative superalgebra with unity ${\mathbf 1}$;
it is the superalgebra
of polynomials in the $a_{\alpha\,i}$ with coefficients
in the group algebra ${\mathbb C}[W(\R)]$ subject to the relations
\bee\label{gh}
\overline g h_\alpha&=&g(h_\alpha)h_\alpha \overline g
%&\ &
\mbox{ for any } g\in W(\R)
                   \mbox{ and } h_\alpha \in {\cal H}^\alpha, \\
\label{rel}
 \!\!\!\!\!\!\!\!\!\!\!\! [ x_\alpha \overline I  , y_\beta \overline I ] &=& \varepsilon_{\alpha\beta}
       \left((\vec x,\, \vec y) {\bar 1}  \overline I +
       \sum_{\vv\in\R} \eta_\vv
\frac {(\vec x,\,\vv)(\vec y,\,\vv)}{(\vv,\,\vv)}{\bar 1} \overline {R_\vv} \right)
% &\ &
\mbox{ for any  $x_\alpha \in {\cal H}^\alpha$ and $y_\beta\, \in {\cal H}^\beta$},
\eee
where $\varepsilon_{\alpha\beta}$
is the antisymmetric tensor, $\varepsilon_{01}=1$, and $\bar 1$ is the unity in ${\mathbb C}[a_{\alpha\,i}]$.
The element   ${\bf 1}=\bar 1 \cdot \overline I$ is the unity  of $H_{W(\R)}(\eta)$.%
\footnote{Clearly, $H_{W(\R)}$ contains  neither $\bar 1\in {\mathbb C}$ nor $\overline I$.}
The action of any operator $g\in End (V)$ is given by a matrix $(g_i^j)$:
\bee\label{ga}
g(a_{\alpha \,i} h^i)&=&a_{\alpha \,i} g_i^j h^j, \quad g_1(g_2( h_\alpha))=(g_1 g_2)(h_\alpha)
\ \mbox{ for any } h_\alpha= a_{\alpha \,i} h^i \in {\cal H}^\alpha,\\
g({\bar 1})&=&{\bar 1}.                                      \label{g1}
\eee
The commutation relations (\ref{rel}) suggest
to define the {\it parity} $\pi$ by setting:
\be
\pi (a_{\alpha\,i}\overline g)=1
\ \mbox{ for any }\alpha,\ i \mbox{ and }g\in \G;
\qquad \pi(\bar 1 \overline g)=0 \ \mbox{ for any } g\in \G.
\nonumber
\ee

We say that $H_{W({R})}(\eta)$ is a {\it superalgebra of observables
of the Calogero model based on the root system $\R$}.
\end{definition}
Clearly,
$H_{W({\R_1+\R_2})}(\eta) = H_{W({\R_1})}(\eta) \otimes H_{W({R_2})}(\eta)$.

%\medskip
These algebras $H_{W({R})}(\eta)$ (with parity forgotten) are particular cases of {\it Symplectic Reflection Algebras}
\cite{sra} and are also known as {\it rational Cherednik algebras} (see, e.g., \cite{ch}).

%\medskip

It follows from eqs. (\ref{ga}) and (\ref{gh})
that if $I\in W(\R)\subset End(V)$ is the unity and
$J=-I\in End(V)$ belongs to $W(\R)$, then
$K:= \bar 1 \overline J\in H_{W(\R)}(\eta)$
is a Klein operator  in $H_{W(\R)}(\eta)$.

%%%%%%%%%%%%%%%%%%%%%%%%%%%%%%%%%%%%%%%%%%%%%%%%%%%%%%%%%%%%%%%%%

\section{Traces and supertraces on $H_{W({R})}(\eta)$}\label{trace}

The following facts were proved in \cite{+1}:

\begin{theorem}\label{th3}
{\it Let the Coxeter group $W(\R)\subset End({\mathbb  R}^N)$
generated by the finite root system $\R\subset{\mathbb  R}^N$ have $
T(\R)$ conjugacy classes without eigenvalue $1$ and $ S(\R)$
conjugacy classes without eigenvalue $-1$.

Then the superalgebra
 $H_{W(\R)}(\eta)$ possesses
$ T(\R)$ independent traces and $ S(\R)$ independent supertraces.
}
\end{theorem}

\begin{theorem}\label{even}
{\it
Each trace and each supertrace on the superalgebta $H_{W(\R)}(\eta)$ is even.
}
\end{theorem}

Theorem \ref{th3} helps to find the numbers $ T(\R)$ and $ S(\R)$ for an arbitrary root
system $\R$.

Theorem \ref{even} implies, clearly, the following statement
\begin{theorem}\label{th4} {\it
In the terms of Theorem \ref{th3}, the following relations are satisfied:}
\bee
 T(\R_1 + \R_2) &=&  T(\R_1)  T(\R_2),\\
 S(\R_1 + \R_2) &=&  S(\R_1)  S(\R_2).
%\nonumber
\eee
\end{theorem}
Therefore, the problem of finding $ T(\R)$ and $ S(\R)$ is reduced to the problem
of finding $ T(\R)$ and $ S(\R)$ for irreducible root systems $\R$.

%The numbers $ S(\R)$ of supertraces for all irreducible root systems
%were preprinted in \cite{stek} and
%the numbers $ T(\R)$ of traces for all irreducible root systems
%were preprinted in \cite{klein}.

Here, the number $ T(\R)$ of traces and
the number $ S(\R)$ of supertraces
for all irreducible root systems are found
and compared. The result is presented in
Sections \ref{=} and \ref{ne=}.

It follows from the results presented in Section \ref{=}, that if
$T(\R)= S(\R)$ for some irreducible root system $\R$, then $-I\in
W(\R)$, and so $H_{W(\R)}(\eta)$ has a Klein operator.

%\newpage

\section{The numbers $ T(\R)$ of traces and $ S(\R)$ of supertraces for irreducible root system $\R$
if $ T(\R)= S(\R)$
}\label{=}

%\vskip -1mm
{\footnotesize
\begin{equation*}
%
%\begin{table}[ht1]
%\tbl{}
{
{\renewcommand{\arraystretch}{0}%
\centering
\begin{tabular}{|c|c|c|c|}
  \hline
\phantom {\small uu} &\phantom {\small uu} &\phantom {\small uu}&\phantom {\small uu}\\
  $\R$ & $ T(\R)= S(\R)$ &
  \begin{tabular}{l}
  presence of $-I$ \\
  in $W(\R)$ proved in
         \end{tabular}
  &
  proof
 in:
%  , see
  \\
\phantom {\small uu} &\phantom {\tiny uu} &\phantom {\tiny uu}&\phantom {\tiny uu}\\
  \hline
\rule{0pt}{2pt}&&\\
  \hline
\phantom {\small uu} &\phantom {\small uu} &\phantom {\small uu}&\phantom {\small uu}\\
  $A_{1}$ &
1 &
\cite{NB}, Table I (XI)   & Appendix \ref{An} \\
\phantom {\small uu} &\phantom {\tiny uu} &\phantom {\tiny uu}&\phantom {\tiny uu}\\
   \hline
\phantom {\small uu} &\phantom {\small uu} &\phantom {\small uu}&\phantom {\small uu}\\
      \begin{tabular}{c}
  $B_{n},C_{n}$%
\\
  \phantom{\tiny uu}%
        \end{tabular}
&
  \begin{tabular}{l}
  the number of partitions  \phantom{f}\\
  of $n$ into the sum of \\
\phantom{\tiny uu}\\
  positive integers
         \end{tabular}
& \cite{NB}, Tables II, III (XI)
& Appendix \ref{Bn} \\
\phantom {\small uu} &\phantom {\tiny uu} &\phantom {\tiny uu}&\phantom {\tiny uu}\\
   \hline
\phantom {\small uu} &\phantom {\small uu} &\phantom {\small uu}&\phantom {\small uu}\\
  $D_{2n}$
&
           \begin{tabular}{l}
   the number of partitions \\
   of $2n$ into the sum of \\
\phantom{\tiny uu}\\
   positive integers with  \\
an   even number of summands
          \end{tabular}
& \cite{NB}, Table IV (XI)
  & Appendix \ref{Dn} \\
\phantom {\small uu} &\phantom {\tiny uu} &\phantom {\tiny uu}&\phantom {\tiny uu}\\
   \hline
\phantom {\small uu} &\phantom {\small uu} &\phantom {\small uu}&\phantom {\small uu}\\
  $E_7$ & 12 &
  \cite{NB}, Table VI (XI)
  & Appendix \ref{7} \\
\phantom {\small uu} &  \phantom {\tiny uu} &\phantom {\tiny uu}&\phantom {\tiny uu}\\
   \hline
\phantom {\small uu} &\phantom {\small uu} &\phantom {\small uu}&\phantom {\small uu}\\
$E_8$ & 30 & \cite{NB}, Table VII (XI)
& Appendix \ref{8} \\
\phantom {\small uu} &\phantom {\tiny uu} &\phantom {\tiny uu}&\phantom {\tiny uu}\\
   \hline
\phantom {\small uu} &\phantom {\small uu} &\phantom {\small uu}&\phantom {\small uu}\\
  $F_4$ & 9 & \cite{NB}, Table VIII (XI)
   & Appendix \ref{4} \\
\phantom {\small uu} &\phantom {\tiny uu} &\phantom {\tiny uu}&\phantom {\tiny uu}\\
   \hline
\phantom {\small uu} &\phantom {\small uu} &\phantom {\small uu}&\phantom {\small uu}\\
  $G_2$ & 3 & \cite{NB}, Table IX (XI)
   & Appendix \ref{2} \\
\phantom {\small uu} &\phantom {\tiny uu} &\phantom {\tiny uu}&\phantom {\tiny uu}\\
   \hline
\phantom {\small uu} &\phantom {\small uu} &\phantom {\small uu}&\phantom {\small uu}\\
  $H_3$ & 4 &  \cite{h4}, p.160; Appendix \ref{3}
   & Appendix \ref{3} \\
\phantom {\small uu} &\phantom {\tiny uu} &\phantom {\tiny uu}&\phantom {\tiny uu}\\
   \hline
\phantom {\small uu} &\phantom {\small uu} &\phantom {\small uu}&\phantom {\small uu}\\
  $H_4$ & 20 & \cite{h4}, Table 3, $K_2$
   & Appendix \ref{H4} \\
\phantom {\small uu} &\phantom {\tiny uu} &\phantom {\tiny uu}&\phantom {\tiny uu}\\
   \hline
\phantom {\small uu} &\phantom {\small uu} &\phantom {\small uu}&\phantom {\small uu}\\
  $I_2(2n)$
&
$ n  $
&
  Appendix \ref{I2}
& Appendix \ref{I2}
\\
\phantom {\small uu} &\phantom {\ uu} &\phantom {\small uu}&\phantom {\tiny uu}\\
   \hline
\end{tabular}%
}
}
%\end{table}
\end{equation*}
}
%%%%%%%%%%%%%%%%%%%%%%%%%%%%%%%%%%%%%%%%%%%%
%\bigskip
%\medskip
\section{The numbers $ T(\R)$ of traces and $ S(\R)$ of supertraces for irreducible root system $\R$
if $ W(R) \not\ni -I $
}\label{ne=}

%\begin{table}[ht2]
%\tbl{}
%{The numbers $ T(\R)$ of traces and $ S(\R)$ of supertraces for irreducible root system $\R$
%if $ W(R) \not\ni -I $.}
{\footnotesize
\begin{equation*}
{
{\renewcommand{\arraystretch}{0}%
\begin{tabular}{|c|c|c|c| }
  \hline
\phantom {\small uu} &\phantom {\small uu} &\phantom {\small uu}&\phantom {\small uu}\\
  $\R$ & $ T(\R)$ & $S(\R)$ &
  proof
 in:
  \\
\phantom {\small uu} &\phantom {\tiny uu} &\phantom {\tiny uu}&\phantom {\tiny uu}\\
  \hline
\rule{0pt}{2pt}&&\\
  \hline
\phantom {\small uu} &\phantom {\small uu} &\phantom {\small uu}&\phantom {\small uu}\\
  $A_0$&
0 & 1
  & Appendix \ref{a0} \\
\phantom {\small uu} &\phantom {\tiny uu} &\phantom {\tiny uu}&\phantom {\tiny uu}\\
   \hline
\phantom {\small uu} &\phantom {\small uu} &\phantom {\small uu}&\phantom {\small uu}\\
  $A_{n-1}$, $n\ge 3$ &
1 &
         \begin{tabular}{l}
  the number of partitions \phantom{f}\\
  of $n$ into the sum of odd \\
\phantom{\tiny uu}\\
  positive integers
   \phantom{f}
         \end{tabular}
  & Appendix \ref{An} \\
\phantom {\small uu} &\phantom {\tiny uu} &\phantom {\tiny uu}&\phantom {\tiny uu}\\
   \hline
\phantom {\small uu} &\phantom {\small uu} &\phantom {\small uu}&\phantom {\small uu}\\
  $D_{2n+1}$
&
           \begin{tabular}{l}
   the number of partitions \\
   of $2n+1$ into the sum of \\
\phantom{\tiny uu}\\
   positive integers with  \\
 an  even number of summands
          \end{tabular}
&
           \begin{tabular}{l}
   the number of partitions \\
   of $2n+1$ into the sum of \\
\phantom{\tiny uu}\\
   positive integers with  \\
an   odd number of summands
          \end{tabular}
  & Appendix \ref{Dn} \\
\phantom {\tiny uu} &\phantom {\tiny uu} &\phantom {\tiny uu}&\phantom {\tiny uu}\\
   \hline
\phantom {\small uu} &\phantom {\small uu} &\phantom {\small uu}&\phantom {\small uu}\\
  $E_6$ & 5  & 9 & Appendix \ref{6} \\
 \phantom {\tiny uu} & \phantom {\tiny uu} &\phantom {\tiny uu}&\phantom {\tiny uu}\\
   \hline
\phantom {\small uu} &\phantom {\small uu} &\phantom {\small uu}&\phantom {\small uu}\\
  $I_2(2n+1)$
& $ n $ & $ n+1 $ &
Appendix \ref{I2} \\
\phantom {\small uu} &\phantom {\ uu} &\phantom {\small uu}&\phantom {\ uu}\\
   \hline
\end{tabular}%
} }
\end{equation*}
}

%\end{table}
%
\vskip 3mm
\indent
{\bf The Weyl superalgebra.}
Let $W_n$ be Weyl superalgebra with $n$ pairs of generating elements:
$W_n={\mathbb C}[a_{\alpha\,i}]$, where $\alpha=0,1$ and $i=1,...,n$,
subject to
relations $[a_{\alpha\,i},\,a_{\beta\,j}]=\varepsilon_{\alpha \beta}\delta_{ij}$
and with parity defined by $\pi(a_{\alpha\,i})=1$.
Clearly, $W_n=(W_1)^{\otimes n}$.
Further, superalgebra $H_{W(A_0)}(\eta)$ does not depend on $\eta$
and since $A_0$ is irreducible, $H_{W(A_0)}(\eta)=W_1$.
So, due to first row of Table \ref{ne=} and Theorem \ref{th4}, the
Weyl superalgebra $W_n$ has 1 supertrace and 0 traces.

%%%%%%%%%%%%%%%%%%%%%%%%%%%%%%%%%%%%%%%%%%%%
%%%%%%%%%%%%%%%%%%%%%%%%%%%%%%%%%%%%%%%%%%%%

\section{Inequality Theorem}

\begin{theorem}\label{th_kl}
{\it
Let $H_{W(\R)}(\eta)$ has $ T(\R)$ traces and $ S(\R)$ supertraces.
Then

i) $ S(\R)>0$,

ii) $ T(\R)\le  S(\R)$,

iii) $ T(\R) =  S(\R)$ if and only if $W(\R)$ contains $-I$.
Equivalently, $ T(\R)= S(\R)$ if and only if $H_{W(\R)}(\eta)$ contains a Klein operator.
}
\end{theorem}

%\medskip
\begin{proof}
%\begin{enumerate}
%\item
Since each group contains unity $I$ and spectrum of $I$ does not contain $-1$,
it follows that $ S(\R)>0$.

%\item
Let $K\in H_{W(\R}(\eta)$ be a Klein operator.
Then $K$ establishes one-to-one correspondence between traces
and supertraces:
\be
\trr(f)=\strr(fK),\ \strr(g)=\trr(gK) \nonumber
\ee

%\item
Let $ T(\R) =  S(\R)$. Then the decomposition of $\R$ in
the sum of irreducible root systems does
not contain root systems from Table \ref{ne=}, namely
$A_0$, $A_n$ for $n\ge 2$, $D_{2n+1}$ for $n\ge 1$ and $E_6$,
because $ T(\R_i) <  S(\R_i)$ for all these root systems.
So, this decomposition contains the root systems listed in Table \ref{=} only,
each of these groups has the element $-I$ and the direct product of
all these $-I$s is $-I$ in $W(\R)$.
%\end{enumerate}

It remains to prove the inequalities
\bee
&& T(A_n)<S(A_n)\qquad\ \ \  \text{ for } n=0 \text{ and } n\ge 2, \label{A}     \\
&& T(D_{2n+1})<S(D_{2n+1}) \ \text{ for }n\ge 1, \label{D}\\
&& T(E_6)<S(E_6). \label{E}
\eee
Inequalities (\ref{A}) and (\ref{E}) manifestly follow from Table \ref{ne=},
and inequality (\ref{D}) follows from Table \ref{ne=} and Lemma \ref{comb}, ii) below.
\end{proof}

\begin{lemma}\!\!\!\!\!\!\!
\footnote{This Lemma is a simple exercise from Partitions Theory,
see, e.g.,  \cite{hall}, \cite{lan} and references therein.
}
\label{comb}%
{\it \hskip 2mm
Let $E(n)$ be the number of partitions of $n$ into the sum of positive integers
with an even number of summands.
Let $O(n)$ be the number of partitions of $n$ into the sum of positive integers
with an odd number of summands.
%\noindent
Then

i)  $E(2k)>O(2k)$ for $k\ge 2$,

ii) $E(2k-1)<O(2k-1)$ for $k\ge 1$,

iii) $E(2)=O(2)$,

iv) $|E(n)-O(n)|=R(n)$, where $R(n)$ is
the number of partitions of $n$ into the sum of different positive odd integers.
}

\end{lemma}
\vskip 3mm
{\begin{proof}
Let $a_{mn}$ be the number of partitions of $n$ into the sum of positive integers
with $m$ summands, $a_{m0}=\delta_{m0}$.
Introduce the generating function%
 $$F(t,x)=\sum_{m,n=0}^\infty a_{mn}t^m x^n.$$
Then
$$\sum_n E(n)x^n = \frac 1 2 (F(t,x)+F(-t,x))|_{t=1}
\mbox{ and }
\sum_n O(n)x^n = \frac 1 2 (F(t,x)-F(-t,x))|_{t=1}.$$

Hence, \[\sum_n (O(n)-E(n))x^n = - F(-t,x)|_{t=1}.\]

Further,
\bee
F(t,x)&=&(1+tx+(tx)^2+(tx)^3+\cdots)(1+tx^2+(tx^2)^2+(tx^2)^3+\dots)
%(1+tx^3+(tx^3)^2+(tx^3)^3+\dots)
\dots \nn
&=&\frac 1 {(1-tx)(1-tx^2)(1-tx^3)\dots} \nonumber
\eee
So
\be\label{frr}
- F(-t,x)|_{t=1}=-\frac 1 {(1+x)(1+x^2)(1+x^3)\dots}
\ee
Multiplying both terms of fraction (\ref{frr}) by
$(1-x)(1-x^2)(1-x^3)\dots$ we obtain
\be\label{fff}
- F(-t,x)|_{t=1}=-\frac {(1-x)(1-x^2)(1-x^3)\dots} {(1-x^2)(1-x^4)(1-x^6)\dots}=-(1-x)(1-x^3)(1-x^5)\dots
\ee
Now, it suffices to notice that eq. (\ref{fff}) can be rewritten in the form
\be%\label{ffff}
\sum_n (O(n)-E(n))x^n=- F(-t,x)|_{t=1}=\sum_{\mbox{$n$ is odd}}R(n)x^n - \sum_{\mbox{$n$ is even}}R(n)x^n \,.
\nonumber
\ee
\end{proof}
}

%%%%%%%%%%%%%%%%%%%%%%%%%%%%%%%%%%%%%%

\section* {Acknowledgments}%

This work was supported
by the Russian Foundation for Basic Research, Grant No.~11-02-00685.

\appendix.

\section{Computing the number of traces and supertraces for all irreducible root systems}

%%%%%%%%%%%%%%%%%%%%%%%%%%%%%%%%%%%%%%%%%%%%%%%%%%%%%%
\subsection{Root system $A_0$}\label{a0}

The Weyl algebra ${\mathbb C}[a,\, a^+]$ generated by elements $a$ and $a^+$ satisfying
the relation $[a,a^+]=1$ may be considered as the algebra of observables of the
Calogero model based on the empty irreducible root system $A_0$, which generates the trivial
group consisting only of the unity $1$. This group has 1 conjugacy class without -1 in its spectrum and 0
conjugacy classes without 1 in its spectrum.
So, this Weyl algebra has 0 traces and 1 supertrace.

\subsection{Root systems $A_{n-1}$ for $n>1$}\label{An}

It is well-known that $W(A_{n-1})=S_n$ and
$V=span(e_1-e_2,\, e_2-e_3,\,...\,,e_{n-1}-e_n)$.
Each element of $S_n$ can be decomposed in the product of cycles
of the form
$$\sigma: e_{i_1}\to e_{i_2}\to \ldots  \to e_{i_k}\to e_{i_1}. $$

\subsubsection{The number of traces}

If $k<n$, then each cycle $\sigma: e_{i_1}\to e_{i_2}\to \ldots  \to e_{i_k}\to e_{i_1} $
has eigenvalue $+1$ with eigenvector
$$e_{i_1}+ e_{i_2}+ \ldots  + e_{i_k} -\frac k n \sum_{s=1}^n e_{s}.$$

The only conjugacy class without eigenvalue $+1$ is the one containing
the cycle of maximal length $n$:
$$\sigma: \ e_{1}\to e_{2}\to \ldots  \to e_{n}\to e_{1} $$
%because its eigenvector $\sum_{i=1}^n e_i$ does not belong to the space $V$.
because its characteristic polynomial has the form
$$
f(t)=1+t+...+t^{n-1}\,.
$$

\subsubsection{The number of supertraces (see also \cite{KV})}

The cycle $\sigma:\  e_{i_1}\to e_{i_2}\to \ldots  \to e_{i_k}\to e_{i_1} $
has eigenvalue $-1$ if and only if $k$ is even. The corresponding eigenvector
has the form
$$
e_{i_1} - e_{i_2}+ \ldots  - e_{i_k}.
$$

So, the number of conjugacy classes without eigenvalue $-1$ is equal to the number of partitions
of $n$ into the sum of positive odd integers \cite{KV}.

\subsubsection{Presence of $-I$ in the group $W(A_{n-1})$, where $n>1$.}

If $A_{n-1}\ni -I$, then $n=2$ (see \cite{NB}, Table I (XI)). The group $W(A_1)$ consists of
two elements: $I$ and $-I$.

\subsection{Root systems $B_n$ and $C_n$}\label{Bn}

The Coxeter group $G=W(B_n)=W(C_n)$ is generated by the permutation
group $S_n$ and
reflections $R_i: \ e_i\to -e_i, \ e_j \to e_j $ for $i \ne j$, see \cite{NB}.
Each element
 $g \in G$ can be represented in the form
 \be\nonumber
 g=\sigma \prod_{i=1}^n R_i ^{\alpha_i},
\mbox{ where $\sigma\in S_n$ and $\alpha_i \in \{0,\, 1 \}$.}
 \ee

 The set $(\sigma,\,\alpha_1,\,...\,\alpha_n)$ unambiguously defines every element of $G$.

 Since each permutation can be decomposed in the product of commuting cycles,
 \be\nonumber
 \sigma=\prod \hat \sigma_k,
 \mbox{~~ where }
\hat \sigma_k: \ e_{i_1} \to e_{i_2} \to \ldots \to e_{i_k} \to e_{i_1},
 \ee
we can introduce what we call {\it $R$-cycles} by the formula
\be\nonumber
\tilde \sigma_k=\hat \sigma_k  R_{i_1}^{\alpha_{i_1}} R_{i_2}^{\alpha_{i_2}} \ldots R_{i_k}^{\alpha_{i_k}}
\ee

So, each element $g\in G$ has the form
\be\label{form}
g=\prod_p \tilde \sigma_p.
\ee

We say that the value $\varepsilon_R(\tilde \sigma)=|\alpha_{i_1} +{\alpha_{i_2}} \ldots +{\alpha_{i_k}}|_{\rm mod\,2 }$
is the {\it $R$-parity of the $R$-cycle $\tilde \sigma$}. Let $l(\tilde \sigma)=k$, where $k$
is the length of the cycle
$\hat\sigma$.

It is easy to prove that an $R$-cycle $\tilde\sigma_1$ is conjugated to an $R$-cycle $\tilde\sigma_2$
if and only if  $l(\tilde\sigma_1)=l(\tilde\sigma_2)$ and $\varepsilon_R(\tilde\sigma_1)=
\varepsilon_R(\tilde\sigma_2)$.

So, a conjugacy class in $G$ is characterized by the numbers $p_1,\,p_2,\,\ldots\,,\,p_n$
and $q_1,\,q_2,\,\ldots\,,\,q_n$, where $p_i$ is the number of $R$-cycles of length $i$ and $R$-parity $0$,
and $q_i$ is the number of $R$-cycles of length $i$ and $R$-parity $1$, in the presentation of $g$ in the form
(\ref{form}).

The numbers $p_i$ and $q_i$ satisfy the relation
\be\label{bbb}
\sum_{i=1}^n (i\,p_i+i\,q_i)=n.
\ee

\subsubsection{The number of traces}

The characteristic polynomial of the $R$-cycle $\tilde\sigma$ has the form
\be\label{charB}
(-1)^{l(\tilde\sigma)}
(t^{l(\tilde\sigma)}-(-1)^{\varepsilon(\tilde\sigma)}).
\ee
It has no root $+1$ if $\varepsilon(\tilde\sigma)=1$.

So, if given conjugacy class has no eigenvalue $+1$, then $p_i=0$ and $\sum_i iq_i =n$.

\subsubsection{The number of supertraces (see also \cite{root})}

The characteristic polynomial of an $R$-cycle $\tilde\sigma$
(\ref {charB}) has no root $-1$ if either $l(\tilde\sigma)$ is even and $\varepsilon(\tilde\sigma)=1$
or if $l(\tilde\sigma)$ is odd and $\varepsilon(\tilde\sigma)=0$.

So, if a given conjugacy class has no eigenvalue $-1$, then
$p_{2k}=0$ and $q_{2k+1}=0$ and eq. (\ref{bbb}) gives
$p_1+2q_2+3p_3+4q_4+...=n$.

\subsubsection{Presence of $-I$ in the group $W(B_{n})$.}
It is easy to see that $-I=\Pi_{i=1}^n R_i$.

\subsection{Root systems $D_n$}\label{Dn}

The Coxeter group $W(D_n)$ is a subgroup of $W(B_n)$, namely,
$g=\prod \tilde\sigma_s$ belongs to $W(D_n)$ if $\left(\sum_s \varepsilon_R (\tilde \sigma_s)
\right)|_{\rm mod \,2}=0$ \cite{NB}.

\subsubsection{The number of traces}

So, $g$ has no eigenvalue $+1$ if $p_i=0$, $\sum_i i\,m_i =n$
and
$\left( \sum_i m_i  \right)|_{\rm mod\, 2}=0$.

This implies that $ T(D_n)$ is equal to the number of partition of $n$ into the sum of
positive integers with an even number of summands.

\subsubsection{The number of supertraces (see also \cite{root})}

Analogously, $ S(D_n)$ is equal to the number of partitions of $n$
   into the sum of positive integers
   with an even number of  even integers.

Clearly, if $n$ is even, then $ S(D_n)$ is equal to the number of partitions of $n$
   into the sum of positive integers
   with an even number of summands,
and if $n$ is odd, then
$ S(D_n)$ is equal to the number of partitions of $n$
   into the sum of positive integers
   with an  odd number of summands.

\subsubsection{Presence of $-I$ in the group $W(D_{n})$.}

If $n$ is even, then $-I=\Pi_{k=1}^n R_k$.

\subsection{Root system $E_6$}\label{6}

The conjugacy classes of the Weyl group $E_6$
are described in Table 9 of \cite{carter}%
\footnote{\label{foot} In \cite{carter}, the Carter diagrams are used
to describe the conjugacy classes of the finite Weyl groups.
Additional results on the Carter diagrams are recently obtained in \cite{carst}.
In particular, two problems are discussed:

i) two Carter diagrams can correspond to
one conjugacy class (\cite{carst}, Theorem 4.1) and

ii)
certain Carter diagrams can correspond to two distinct conjugacy classes
(\cite{carst}, Theorem 6.5).
}.

\subsubsection{The number of traces}

The following 5 classes do not have the root $+1$:
$$
     A_2^3,
\   A_5 \times A_1,
  \   E_6,
\   E_6(a_1),
\   E_6(a_2)%
.
$$

\subsubsection{The number of traces}

The following 9 classes do not have the root $-1$:
$$
   \phi,
\    A_2,
\   A_2^2,
  \   A_4,
\   D_4(a_1),
\   A_2^3,
\   E_6,
\   E_6(a_1),
\   E_6(a_2).
$$

\subsection{Root system $E_7$}\label{7}

The conjugacy classes of the Weyl group $E_7$
are described in Table 10 of \cite{carter}${}^{\ref{foot}}$.%

\subsubsection{The number of traces}

The last 12 classes in Table 10 of \cite{carter} do not have the root $+1$:
\begin{eqnarray*}
&&
\    A_1^7,
\   A_3^2 \times A_1,
  \   A_5 \times A_2,
\   A_7,
\   D_4 \times A_1^3,
\  D_6 \times A_1,
\  D_6(a_2) \times A_1,
\\&&
\  E_7,
\   E_7(a_1),
\   E_7(a_2),
\   E_7(a_3),
\   E_7(a_4).
\end{eqnarray*}

\subsubsection{The number of supertraces}

The following 12 classes do not have the root $-1$:
$$
   \phi,
\    A_2,
\   A_2^2,
  \   A_4,
\   D_4(a_1),
\   A_2^3,
\  A_4 \times A_2,
\  A_6,
\   D_6(a_1),
\   E_6,
\   E_6(a_1),
\   E_6(a_2).
$$

\subsection{Root system $E_8$}\label{8}

The conjugacy classes of the Weyl group $E_8$ are described in Table 11 of \cite{carter}%
${}^{\ref{foot}}$.%

\subsubsection{The number of traces}

The last 30 classes in Table 11 of \cite{carter} do not have the root $+1$:
\begin{eqnarray*}
&&
    A_1^8,
\   A_2^4,
\   A_3^2 \times A_1^2,
\   A_4^2,
\   A_5 \times A_2 \times A_1 ,
\   A_7 \times A_1,
\   A_8,
\   D_4 \times A_1^4,
\   D_4^2,
\   D_4(a_1)^2,
\\
&&
    D_5(a_1) \times A_3,
\   D_6 \times A_1^2,
\   D_8,
\   D_8(a_1),
\   D_8(a_2),
\   D_8(a_3),
\   E_6 \times A_2,
\   E_6(a_2) \times A_2,
\\
&&
    E_7 \times A_1,
\   E_7(a_2) \times A_1,
\   E_7(a_4) \times A_1,
\   E_8,\  E_8(a_i)\  (i = 1,\,...\,8).
\end{eqnarray*}

\subsubsection{The number of supertraces}

The following 30 classes do not have the root $-1$:
\begin{eqnarray*}
&&
\phi,
\    A_2,
\   A_2^2,
\   A_4,
\   D_4(a_1),
\   A_2^3,
\  A_4 \times A_2,
\  A_6,
\  D_4(a_1) \times A_2,
\   D_6(a_1),
\\
&&
E_6,
\   E_6(a_1),
\   E_6(a_2),
\   A_2^4,
\   A_4^2,
\   A_8,
\   D_4(a_1)^2,
\   D_8(a_1),
\   D_8(a_3),
\\
&&
E_6 \times A_2,
\   E_6(a_2) \times A_2,
\   E_8,\  E_8(a_i)\  (i = 1,\,...\,8).
\end{eqnarray*}

\subsection{Root system $F_4$}\label{4}

The conjugacy classes of the Weyl group $F_4$
are described in Table 8 of \cite{carter}${}^{\ref{foot}}$.%

\subsubsection{The number of traces}

The last 9 classes in Table 8 of \cite{carter} do not have the root $+1$:
$$
  A_1^4,
\ A_2 \times \tilde A_2,
\ A_3 \times \tilde A_1,
\ C_3 \times A_1,
\ D_4,
\ D_4(a_1),
\ B_4,
\ F_4,
\ F_4(a_1).
$$

\subsubsection{The number of supertraces}

The following 9 classes do not have the root $-1$:
$$\phi,\ A_2, \ \tilde A_2,\ B_2,\ A_2 \times \tilde A_2,\ D_4(a_1),
\ B_4,\ F_4,\ F_4(a_1).$$

\subsection{Root system $G_2$}\label{2}

The conjugacy classes of the Weyl group $G_2$
are described in Table 7 of \cite{carter}${}^{\ref{foot}}$.%

\subsubsection{The number of traces}

The last 3 classes in Table 7 of \cite{carter} do not have the root $+1$:
$$
  A_1 \times \tilde A_1,
\ A_2,
\ G_2,.
$$

\subsubsection{The number of supertraces (see \cite{i2})}

The following 3 classes in Table 7 of \cite{carter} do not have the root $-1$:
$$
 \phi,
\ A_2,
\ G_2,.
$$

The fact that $G_2$ has 3 conjugacy classes without eigenvalue $+1$ and
3 conjugacy classes without eigenvalue $-1$
can be derived also from Appendix \ref{I2} because
$G_2=I_2(6)$.

\subsection{Root system $H_3$}\label{3}

Let $k=\frac12(\sqrt {5} +1)$.
Then the reflections
\begin{eqnarray*}
a=
\left(
\begin{array}{ccc}
1& 0& 0 \\
0 &-1& 0 \\
0 &0& 1
\end{array}
\right),
\qquad
b=\frac 1 2
\left(
\begin{array}{ccc}
1& k& k-1 \\
k &1-k& -1 \\
k-1 &-1& k
\end{array}
\right),
\qquad
c =
\left(
\begin{array}{ccc}
-1& 0& 0 \\
0 &1& 0 \\
0 &0& 1
\end{array}
\right)
\\
\end{eqnarray*}
corresponding to the roots $\vec e_2$,\
$\frac 1 2 (-\vec e_1 + k\vec e_2 + k^{-1}\vec e_3)$
and $\vec e_1$, respectively,
satisfy the relations
$$a^2 = b^2 = c^2 = 1, \ \ (ab)^5 = (bc)^3 = (ac)^2 = 1$$
and generate the Coxeter group $H_3$.

As
$H_3 = S^e_5 \times C_2$ (See \cite{NB},  Ch.VI, Sect. 4, Ex.11 d),  p. 284;
\cite{h4}, p.160 and references therein),
where $S^e_5$ is the group of even permutations of 5 elements,
$C_2=\{1,-1\}$, it follows that
$H_3$ has 10 conjugacy classes, 5 with a positive determinant
and 5 with a negative one.
(Observe that $|S^e_5| = 60$, hence  $|H_3| = 120$.)

The conjugacy classes with positive determinant are
described by their representatives

\begin{equation}\label{table}
%\centerline{
\begin{tabular}{|l|l|}
   \hline
\phantom {\tiny uu} &\phantom {\tiny uu}\\
The representative & The characteristic polynomial \\
\phantom {\tiny uu} &\phantom {\tiny uu}\\
   \hline
\phantom {\tiny uu} &\phantom {\tiny uu}\\
$ I$ &  $( 1 - t )^3$
\\
$ ac$  & $( 1 - t )( 1 + t )^2$
\\
$ bc$  & $( 1 - t )( t^2 + t + 1 )$
\\
 $ab$  & $( 1 - t ) [ t^2 + ( 1 - k ) t + 1 ]$
\\
$ abab$  & $( 1 - t ) ( t^2 + k t + 1 )$
\\
\phantom {\tiny uu} &\phantom {\tiny uu}\\
\hline
\end{tabular}
%}
\end{equation}%\bigskip

\subsubsection{The number of supertraces}

Only four conjugacy classes have no roots $-1$, their representatives are
$I$, $bc$, $ab$ and $abab$.

\subsubsection{The number of traces}

Each of the characteristic polynomials (\ref{table}) has the root $+1$.
Besides, the conjugacy class with representative $-ac$ has the root $+1$,
so only four conjugacy classes with negative determinant have no roots $+1$.
Their representatives are
$-I$, $-bc$, $-ab$ and $-abab$.

So, the number of conjugacy classes without root $+1$ is equal to 4.

\subsubsection{Presence of $J=-I$ in $W(H_3)$.}

The group $H_3$ contains the element $J=-I=(ababc)^3$ (see \cite{w0}, p.11).

%\newpage

\subsection{Root system $H_4$}\label{H4}

According to \cite{h4}, all 34 conjugacy classes of the Coxeter group  $H_4$ are described by
their representatives acting on the space of quaternions
\bee\label{lr}
g_{lr}:\ \ x \mapsto l x r^*,
\\
g^*_{p}:\ \ x \mapsto p x^*.\label{p}
\eee
All 25 pairs of unit quaternions $l$ and $r$ and 9 unit quaternions $p$
are listed in Table 3 of \cite{h4}.

\subsubsection{The number of traces}

Each operator (\ref{p}) has the root $+1$.
Indeed, the equation
$
px^*=x
$
has a nonzero solution
$
x= 1+p
$
if $p\ne -1$,
and $x$ is an arbitrary imaginary quaternion if $p=-1$.

Each operator (\ref{lr}) does not have the root $+1$ if and only if
\be\label{ne}
l_0-r_0\ne 0.
\ee
Indeed, the determinant of the map
$x\mapsto lx-xr$ is equal to $4(l_0 - r_0)^2$.

There are 20 pairs $(l,r)$ in Table 3 of \cite{h4}
satisfying the condition (\ref{ne}), namely,
$K_i$ with $i$ = 2, 5, 7, 9 to 25.

\subsubsection{The number of supertraces}

Each operator (\ref{p}) has the root $-1$.
Indeed, the equation
$
px^*=-x
$
has nonzero solution
$
x= 1-p
$
if $p\ne 1$,
and $x$ is an arbitrary imaginary quaternion if $p=1$.

Each operator (\ref{lr}) does not have the root $-1$ if and only if
\be\label{ne+}
l_0+r_0\ne 0.
\ee
Indeed, the determinant of the map
$x\mapsto lx-xr$ is equal to $4(l_0 + r_0)^2$.

There are 20 pairs of $l,r$ in Table 3 of \cite{h4}
satisfying the condition (\ref{ne+}), namely,
$K_i$ with $i$ = 1, 4, 6, 8, 10 to 25.

\subsubsection{Presence of element $J=-I$ in $W(H_4)$.}

The element $K_2$ in Table 3 of \cite{h4} with $l=-r=1$ is
$-I$ in $H_4$.

\subsection{Root systems $I_2(n)$}\label{I2}

It is convenient to use ${\mathbb C}$ instead of ${\mathbb  R}^2$ to
describe $W(I_2(n))$. The root system $I_2(n)$ contains $2n$ vectors
$v_k=\exp(\pi i k/n)$, where $k=0,1,...,2n-1$. The corresponding
Coxeter group $W(I_2(n))$ has $2n$ elements, $n$ reflections $R_k$
acting on $z,\ z^* \in{\mathbb C}$ as follows
\bee%\label{Rk}
R_k z&=& -z^* v_k^2 R_k, \nn
R_k z^* &=& -z {v_k^*}^2 R_k,\qquad k\in {\mathbb Z}_n  \nonumber
\eee
and $n$ elements of the form $S_k=R_k R_0$, where $S_0$ is the unity in $W(I_2(n))$.
These elements satisfy the following relations
\be %\label{RR}
R_k R_l = S_{k-l},\qquad
S_k S_l = S_{k+l},\qquad
R_k S_l = R_{k-l},\qquad
S_k R_l = R_{k+l}.
\nonumber
\ee

Obviously, the reflections $R_{2k}$ lie in one conjugacy class
and $R_{2k+1}$ lie in another one if $n$ is even. If $n$ is odd, then all reflections
$R_k$ lie in one conjugacy class. Each reflection has both eigenvalues $+1$ and $-1$.

The rotation $S_k$ has no eigenvalues $-1$ if $k \ne n/2$,
and has no eigenvalues $+1$ if $k \ne 0$.
If $n$ is even, then $S_{n/2}=-I$.

Rotations $S_k$ and $S_{-k}$ form a conjugacy class.

So, the number of conjugacy classes without $+1$ is equal to $\left[\frac n 2 \right ]$, and
the number of conjugacy classes without $-1$ is equal to $\left[\frac {n+1} 2 \right ]$, see  \cite{i2}.

More properties of the algebra $H_{1,\eta} (I_2(2m + 1))
:=
H_{W({I_2(2m + 1)})}(\eta)$ see in \cite{i22}.


\begin{thebibliography}{00}

\bibitem{NB} N.~Bourbaki, {\it Groupes et alg\`ebres de Lie},
                         ch.IV--VI, Hermann, Paris, 1968.

\bibitem{carter} R.W.~Carter,
       \lq\lq Conjugacy classes in the Weyl group", Compositio Mathematica,
       {\bf 25} $\mathsf n^o$1, (1972) 1--59.

\bibitem{sra} P.~Etingof and V.~Ginzburg,
``Symplectic reflection algebras, Calogero--Moser
space, and deformed Harish--Chandra homomorphism", Inv. Math. 147(2002), 243--348;
arXiv:math.AG/0011114.

\bibitem{ch} P.~Etingof, Xiaoguang Ma,
             ``Lecture notes on Cherednik algebras",
             arXiv:1001.0432v4.

\bibitem{w0} W.N.~Franzsen, \lq\lq Automorphisms of Coxeter Groups",
             A thesis for the degree of Doctor of Philosophy,
             School of Mathematics and Statistics University of Sydney,
             January, 2001.

\bibitem{h4} L.C.~Grove,
       \lq\lq The characters of the hecatonicosahedroidal group",
       Journal f\" ur die reine und angewandte Mathematik,
        {\bf 265} (1974) 160--169.

\bibitem{hall} M.~Hall, {\it Combinatorial theory}, 2nd edition,
               Wisley, New York (1986), 464 pp.

\bibitem{JH} J.E.~Humphreys, {\it Reflection groups and Coxeter
             groups}, Cambridge University Press, 1990.

\bibitem{root} S.E.~Konstein,
              Teor.~Mat.~Fiz., {\bf 111} (1997) 252--262.

\bibitem{i2} S.E.~Konstein,
       \lq \lq Rational Calogero models based on rank-2 root systems:
       supertraces on the superalgebras of observables",
       arXiv:math/9801001.

\bibitem{Ko}
S.E.~Konstein, \lq\lq Supertraces on the Superalgebra of Observables of
                Rational Calogero Model based on the Root System",
               arXiv:math-ph/9904032.

\bibitem{stek} S.E.~Konstein and R.~Stekolshchik,
             \lq\lq The Number of Supertraces on the Superalgebra of
              Observables of Rational Calogero Model
              based on the Root System",
              arXiv:0811.2487.

\bibitem{+1}
S.E.~Konstein and I.V.~Tyutin,
       \lq\lq Traces on the Superalgebra of Observables of
       Rational Calogero Model based on the Root System",
       Journal of Nonlinear Mathematical Physics,
       Vol. 20, No. 2 (June 2013), 271-294;
        arXiv:1211.6600.

\bibitem{i22}
S.E.~Konstein and I.V.~Tyutin  \lq\lq Ideals generated by traces or by
supertraces in the symplectic reflection algebra $H_{1,\nu} (I_2(2m + 1))$\rq\rq,
Journal of Nonlinear Mathematical Physics, 24:3, (2017) 405-425,
DOI: 10.1080/14029251.2017.1341702

\bibitem{KV} S.E.~Konstein and M.A.~Vasiliev,
              J.~Math.~Phys. {\bf 37} (1996) 2872

\bibitem{lan} S.K.~Lando, {\it Lectures on generating functions}.
              Series: Student mathematical library, v. 23,
              Published 2003 by American Mathematical Society in Providence, RI,  144 pp.

\bibitem{OP} M.A.~Olshanetsky and A.M.~Perelomov,
              Phys. Rep., {\bf 94} (1983) 313.

\bibitem{sosrus} {\em Seminar on supersymmetry (v.  $1$.
                 Algebra and Calculus: Main chapters)},
                 (J.~Bernstein, D.~Leites, V.~Molotkov, V.~Shander)
                 ed. by D.~Leites, MCCME, Moscow, {\bf 2011}, 410 pp
                 (in Russian; a version in English is in preparation).

\bibitem{raf} R.~Stekolshchik,
              {\it Notes on Coxeter Transformatios and the McKay Correspondence}.
              Series: Springer Monographs in Mathematics 2008, XX, 239 pp.

\bibitem{carst} R.~Stekolshchik,
                  \lq\lq Root systems and diagram calculus.
                  I. Regular extensions of Carter diagrams and the
                  uniqueness of conjugacy classes",
                  arXiv:1005.2769v6.

\bibitem{V} M.A.~Vasiliev, JETP Letters, {\bf 50} (1989) 344--347;
              Int. J. Mod. Phys. {\bf A6} (1991) 1115.

\end{thebibliography}
\end{document}